\documentclass[11pt]{amsart}
\usepackage{amssymb,amscd}
\input{xy}
\xyoption{all}

\newtheorem{lemma}{Lemma}[section]
\newtheorem{corollary}[lemma]{Corollary}
\newtheorem{theorem}{Theorem}
\newtheorem{proposition}[lemma]{Proposition}
\theoremstyle{definition}
\newtheorem{example}{Example}
\newtheorem{remark}{Remark}

\newtheorem{definition}[lemma]{Definition}
\newtheorem*{ack}{Acknowledgement}
\def\set#1{\{#1\}}
\def\abs#1{{\left|#1\right|}}
\def\angl#1{{\left< #1\right>}}
\def\Mo{{\mathcal Z}}
\def\CR{{{\mathcal C}_\rho}}

\def\I{{\mathcal I}}
\def\Z{{\mathbb Z}}

\def\R{{\mathbb R}}

\def\k{{\bf k}}
\def\a{{\mathfrak a}}
\def\wf{{\widetilde{f}}}
\def\Ext{{\rm Ext}}
\def\Hom{{\rm Hom}}
\DeclareMathOperator{\ann}{ann}
\def\Tor{{\rm Tor}}

\def\link{{\rm link}}
\def\supp{{\rm supp}}

\def\vareps{{\varepsilon}}
\def\rhot{{\widetilde{\rho}}}
\def\gr{{\rm Gr}}
\def\aug{{\bf 1}}
\def\bfcdot{{\mathbf\cdot}}
\newdir^{ (}{{}*!/-5pt/@^{(}}
\subjclass[2000]{Primary
20F36,  
Secondary
13F55,  
20J05.  
}
\title[Homology of subgroups of right-angled Artin groups]
{Homology of subgroups of right-angled Artin groups}

\author[G. Denham]{Graham Denham}
\address{Department of Mathematics, University of Western Ontario,
London, ON  N6A 5B7}
\urladdr{{http://www.math.uwo.ca/\~{}gdenham}}
\thanks{Partially supported by a grant from NSERC of Canada.}

\begin{document}
\begin{abstract}
We describe the (co)homology of a certain family of 
normal subgroups of right-angled Artin
groups that contain the commutator subgroup, as modules over the quotient
group.  We do so in terms of (skew) commutative algebra of 
squarefree monomial ideals.  
\end{abstract}
\maketitle
\section{Introduction}
Let $\Gamma$ be a simple graph on $n$ vertices $V=\set{v_1,v_2,\ldots,v_n}$.
The right-angled Artin group (or graph group) $G_\Gamma$ is the group
with generators $V$ and relations $v_iv_j=v_jv_i$ for each edge $v_iv_j$ in
$\Gamma$.  Charney and Davis~\cite{chardav95} showed that such groups admit a 
finite classifying space, a subcomplex of the $n$-torus first introduced
in \cite{kr80}.  It follows
that the cohomology ring of $G_\Gamma$ has an elegant description as
an exterior Stanley-Reisner ring.

Bestvina and Brady~\cite{bb97} make use of a particular subgroup
of $G_\Gamma$ with infinite cyclic quotient as examples to 
distinguish finiteness
properties:
the kernel of the map to $\Z$ sending each $v_i$ to the generator
$1$
is finitely generated if and only if $\Gamma$ is connected and finitely
presented if and only if $\Gamma$ is simply connected.  Moreover,
they show that this subgroup is $FP_k$ if and only
if the flag complex $K_\Gamma$ of $\Gamma$ is homologically $k$-connected.
The cohomology ring of this subgroup is computed in 
\cite{learsaad06,pasu06b} when
it is finite-dimensional, by relating it to the simplicial topology of the 
flag complex $K_\Gamma$.\nocite{pasu06}

More generally, for an integer $m$ let $\rho\colon G_\Gamma\to\Z^m$ be a 
surjective group homomorphism with the property that $\rho(v_i)$, for
each $i$, is a generator of the abelian group $\Z^m$.  We shall call
such a map $\rho$ a {\em coordinate homomorphism}, and denote its
kernel by $N_\rho=N_{\Gamma,\rho}$, 
the {\em coordinate subgroup}.  From Meier, Meinert
and VanWyk's calculation of the Bieri-Neumann-Strebel invariants of
right-angled Artin groups~\cite{mmv98}, it follows that the
homology groups of $N_\rho$ are not finitely generated except under very
restrictive hypotheses.

The point of view of this paper is that, nevertheless, $H_*(N_\rho,\Z)$
is a finitely-generated module over the group ring $\Z[\Z^m]$, and so
is amenable
to description in terms of the graph $\Gamma$ 
via combinatorial commutative algebra.
Accordingly, we compute this module in terms of the exterior Stanley-Reisner
ring of the clique complex of the graph $\Gamma$.  In such terms, for example,
one can determine the Krull dimension of each module $H_p(N_\rho,\Z)$.  

Under the additional hypothesis 
that the complex $K_\Gamma$ is Cohen-Macaulay, these results
can be made more explicit via Bernstein-Gelfand-Gelfand duality 
(Section~\ref{sec:duality}).  In particular, we find that
$H^q(G_\Gamma,
\Z[G_\Gamma^{ab}])$ is zero except for $q=d+1$ if and only if
$K$ is Cohen-Macaulay of dimension $d$.  This is an abelian version 
of a result of Brady and Meier~\cite{BM01}; in this case, we are also
able to describe the dualizing module explicitly.

\section{Classifying spaces}
The construction used here is a generalization of constructions that 
appear independently in the work of various authors.  In the context
of right-angled Artin groups, the idea originates with Charney and
Davis~\cite{chardav95}.  The language of partial product complexes
is convenient,
however; details and further references may be found in \cite{desu06b}.
(These are also known as generalized moment-angled complexes; see 
\cite{Str99,bp02}.)
\subsection{Partial product complexes}
\begin{definition}
\label{def:zkx}
Let $X$ be a space, and $A\subset X$ a non-empty 
subspace.  Given a simplicial complex $K$ on vertex set 
$[n]=\set{1,2,\ldots,n}$, define $\Mo_K(X,A)$ to be 
the following subspace of the cartesian product $X^{\times n}$:
\begin{equation}
\label{eq:zkx}
\Mo_K(X,A)=\bigcup_{\sigma\in K}  (X,A)^{\sigma}, 
\end{equation}
where $(X,A)^{\sigma}= \{ x \in X^{\times n} \mid x_i \in A \text{ if } 
i\notin \sigma\}$.
\end{definition}

If $X$ is a pointed space, let $\Mo_K(X)=\Mo_K(X,*)$.  For example,
$\Mo_K(S^1)$ is a subcomplex of the $n$-torus $(S^1)^{\times n}$.

\subsection{Right-angled Artin groups}
If $\Gamma$ is a graph with vertices $V(\Gamma)=\set{v_1,\ldots,v_n}$ and
edges $E(\Gamma)$, recall
the {\em right-angled Artin group} $G_\Gamma$ is defined by the presentation
\begin{equation}\label{eq:defartin}
G_\Gamma=
\left< v_1,\ldots,v_n\mid v_iv_j=v_jv_i\text{ for each $v_iv_j\in E(\Gamma)$}
\right>
\end{equation}

Also recall that if $\Gamma$ is a graph, its {\em clique complex} $K_\Gamma$
is the simplicial complex with vertices $V(\Gamma)$ and simplices
$\sigma$, for all $\sigma\subseteq V(\Gamma)$ with the property that
each pair of vertices of $\sigma$ are connected by an edge.  If a simplical
complex $K$ is the clique complex of its $1$-skeleton, $K$ is called a 
flag complex.
Then the construction of Definition~\ref{def:zkx} recovers the construction
of \cite{chardav95}:

\begin{proposition}[\cite{chardav95}]\label{prop:pi1}
Let $K$ be a simplicial complex and let $\Gamma=K^{(1)}$, its $1$-skeleton.
Then $G_\Gamma\cong\pi_1(\Mo_K(S^1),*)$.  If, further, $K$ is a flag complex,
then $\Mo_K(S^1)$ is an Eilenberg-Maclane space for $G_\Gamma$.
\end{proposition}

For example, $\pi_1(\Mo_K(S^1))$ is abelian if and only if the $1$-skeleton
of $K$ is a complete graph.  The clique complex of a complete graph is simply
a full simplex $K=\Delta^{n-1}$ on $n$ vertices, in which case $\Mo_K(S^1)$
is the $n$-torus $(S^1)^{\times n}$.

Since the main tool used here is the space $\Mo_K(S^1)$ (rather than the
group $G_\Gamma$), it will be natural to consider the homology of spaces,
rather than groups; as the Proposition indicates, we will recover the
case of groups by specializing to those $K$ which are flag complexes.
We shall correspondingly regard simplicial complexes as our primary objects,
rather than graphs.

\subsection{Coordinate homomorphisms}
Here we identify the coordinate subgroups, defined in the
Introduction, in terms of our topological construction.

Such subgroups are a special case of a more general construction.
If $f\colon K\to L$ is a map of simplicial complexes sending vertices
$[n]$ to vertices $[m]$, there is a natural map $\Mo_f:\Mo_K(S^1)\to
\Mo_L(S^1)$ obtained by restricting a map $\overline{f}\colon(S^1)^{\times n}
\to(S^1)^{\times m}$.  Here, $\overline{f}(x)_j=\prod_{i\colon f(i)=j}x_i$;
see \cite[Lemma~2.2.2]{desu06b} for details.

Using Proposition~\ref{prop:pi1}, it is routine to translate this to
a statement about right-angled Artin groups.
\begin{proposition}\label{prop:inducedhom}
If $f\colon K\to L$ is a map of simplicial complexes, the induced 
map of fundamental groups $\Mo_f^\sharp\colon G_{K^{(1)}}\to G_{L^{(1)}}$
sends the $i$th generator of $G_{K^{(1)}}$ to the $f(i)$th generator in 
$G_{L^{(1)}}$.
\end{proposition}

In particular, if $\Gamma=K^{(1)}$ has $n$ vertices, then 
the abelianization of a right-angled Artin group $G_\Gamma$ is clearly
$\Z^n$, and the abelianization map $G_\Gamma\to\Z^n$ is obtained by
choosing $L=\Delta^{n-1}$, the full simplex on $n$ vertices.
More generally, we consider the following case:
\begin{definition}\label{def:cmap}
Let $K$ be a simplicial complex on $n$ vertices, and
let $f\colon[n]\to[m]$ be a surjective function on sets.  Then 
$f$ extends uniquely 
to a map of simplicial complexes $f\colon K\to\Delta^{m-1}$.
We will call $f$ a coordinate map.
\end{definition}

It follows from Proposition~\ref{prop:inducedhom} that a coordinate map
$f\colon K\to[m]$ induces a coordinate homomorphism $\Mo_f^\sharp$, and
this homomorphism factors through the abelianization of $G_{K^{(1)}}$:
\begin{equation*}
\xymatrix@R12pt@C10pt{
\Mo_f^\sharp\colon G_{K^{(1)}}\ar[rr]
\ar[dr] && \Z^m\\
&\Z^n\ar[ur]_{\overline{f}^\sharp}
}
\end{equation*}
Conversely, by comparing definitions one finds that all 
coordinate homomorphisms arise in this way.

\begin{example}\label{ex:univab}
At the one extreme, we may take $f$ to be the identity map.  Then
the homomorphism $\Mo_f^\sharp\colon G_{K^{(1)}}\to\Z^n$ is the abelianization,
and its kernel $N_f$ is the commutator subgroup of the right-angled
Artin group.
\end{example}
\begin{example}
On the other hand, the (unique) map $f\colon [n]\to[1]$ induces a homomorphism
$\Mo_f^\sharp\colon G_{K^{(1)}}\to\Z$ sending each $v_i$ to $1$.  The kernel
is the Bestvina-Brady group $H_\Gamma$ considered in 
\cite{bb97,mmv98,mmv01,learsaad06,pasu06b}.
\end{example}

\subsection{Abelian covers}
Let $\pi:\R\to S^1$ be the universal cover of $S^1$, sending $\Z$ to
the basepoint $*$.
By \cite[Lemma~2.9]{desu06b}, the map of pairs $\pi:(\R,\Z)\to (S^1,*)$
induces a fibration
\begin{equation}\label{eq:fib}
\Z^n\rightarrow \Mo_K(\R,\Z)\stackrel{\Mo_\pi}\rightarrow \Mo_K(S^1),
\end{equation}
which is in fact the universal cover of $\Mo_K(S^1)$.

By covering space theory, if $f:K\to[m]$ is a coordinate
homomorphism, then $N_f$ is the fundamental group of the
fibred coproduct $\Mo_K(\R,\Z)\times_{\Z^n}\Z^m$,
where $\Z^n$ acts on $\Mo_K(\R,\Z)$ by deck transformations, and
on $\Z^m$ by the induced map $\overline{f}^\sharp$.

\subsection{CW-complexes}\label{ss:cw}
Combinatorially explicit cell structures are available for each
complex.  First, the standard cell structure on the torus $(S^1)^{\times n}$
restricts to give a cell structure for $\Mo_K(S^1)$ (see~\cite{kr80}).
From the definition \eqref{eq:zkx}, 
the cells are naturally labelled by the simplices of $K$.  For each
simplex $I\in K$ with $k$ vertices, let $\vareps_I$ denote the corresponding
$k$-cell.  Note that the complex is minimal in the sense that the attaching
maps are zero.

We need the following notation.  Let $E=\Z[e_1,\ldots,e_n]$ denote the
exterior algebra, which we regard as a 
graded-commutative Hopf algebra with generators in degree $1$.
Let $\set{\vareps_i}$ denote the $\Z$-dual basis to the generators, so that
$E^*=\Z[\vareps_1,\ldots,\vareps_n]$ is also an exterior algebra.

Now let $\Z\angl{K}=E/J_K$, the exterior Stanley-Reisner ring of $K$, where $J_K$ 
is the
ideal generated by monomials indexed by nonfaces of $K$.  Then
$\Z\angl{K}^*$ is a sub-coalgebra of $E^*$, spanned by monomials $\vareps_I$,
where we set $\vareps_I:=\vareps_{i_1}\cdots\vareps_{i_k}$ if
$I=\set{i_1,\ldots,i_k}$ is a simplex of $K$.  We identify
$\Z\angl{K}^*$ with the cellular chain complex of $\Mo_K(S^1)$ (with 
zero differential).  The reason for this notation is the following
foundational result.
\begin{theorem}[\cite{kr80}]\label{th:EJ}
Let $K$ be a simplicial complex.  Then $H^*(\Mo_K(S^1),\Z)\cong \Z\angl{K}$
as graded rings.
\end{theorem}

Second, we require similar cell structure for the complex $\Mo_K(\R,\Z)$,
It turns out that, in the flag-complex setting, the construction below
is really the (abelianized) Salvetti complex
for the right-angled Artin group, as explained and generalized by Charney
and Davis~\cite{chardav95}.  

To start, label the zero
cells of $\R$ by $\set{x^i:i\in\Z}$, and the $1$-cells by
$\set{\vareps\cdot x^i:i\in\Z}$ so that $\partial(\vareps\cdot x^i)=
x^{i+1}-x^i$
for each $i$.  Extending this to the product structure on $\R^n$
and restricting to the subcomplex $\Mo_K(\R,\Z)$ gives a complex
with $k$-cells labelled in a natural way by 
$$
\set{\vareps_I\cdot x^\alpha\colon I\in K,\abs{I}=k,\hbox{~and~}
\alpha\in\Z^n},
$$
where we regard $x^\alpha$ equivalently as an element of $\Z^n$
(written multiplicatively) or a Laurent monomial $x^\alpha:=x_1^{\alpha_1}
\cdots x_n^{\alpha_n}$.
It is straightforward to check that the differential is given by 
\begin{equation}\label{eq:diff0}
\partial(\vareps_Ix^\alpha)=\sum_{i\in I}(-1)^{\sigma(I,i)}
\vareps_{I-\set{i}}\cdot x^{\alpha+\chi_i}-\vareps _{I-\set{i}}\cdot x^\alpha,
\end{equation}
where $\sigma(I,i)$ is the number of elements preceding $i$ in $I$ 
(in the standard order) and $\chi_i$ is the $i$th coordinate vector
in $\Z^n$.  Make the identification
\begin{equation}\label{eq:cell}
C_*^{\rm cell}(\Mo_K(\R,\Z))\cong \Z\angl{K}^*\otimes_\Z \Z[\Z^n];
\end{equation}
then some useful properties follow immediately from the construction.
(Recall $v_i$ denotes the $i$th standard generator of $\pi_1(\Mo_K(S^1))$
from \eqref{eq:defartin}.)
\begin{proposition}\label{prop:equivariant}
For any $K$, 
\begin{enumerate}
\item The universal cover $\Mo_\pi$ of \eqref{eq:fib} cellular and
satisfies $\Mo_\pi(\vareps_I x^\alpha)=\vareps_I$ for all choices of 
$I$ and $\alpha$.
\item The cell structure \eqref{eq:cell} 
is $\Z^n$-equivariant, and $v_i$ acts by
multiplication by $x_i$ for $1\leq i\leq n$;
\item The differential \eqref{eq:diff0} is induced by 
$\partial(\vareps_i)=x_i-1$ for each $1\leq i\leq n$, together with
$\Z[\Z^n]$-linearity and the Leibniz rule on $\Z\angl{K}^*$;
\end{enumerate}
\end{proposition}
\section{(Co)homology of abelian covers}
In this section, we describe the homology of the coordinate subgroups
in terms of commutative algebra.  If $f$ is a coordinate map 
(Def.~\ref{def:cmap}), then $\Z[\Z^m]$ is a
$\pi_1(\Mo_K(S^1))$-module via the homomorphism $\Mo_f^\sharp$.  Then
$$
H_*(\Mo_K(S^1)\times_{\Z^n}\Z^m,\Z) \cong
H_*(\Mo_K(S^1),\Z[\Z^m]);
$$
in the case where $K$ is a flag complex and $\Gamma=K^{(1)}$, this is
simply Shapiro's Lemma:
$$
H_*(N_f,\Z)\cong H_*(G_\Gamma,\Z[\Z^m]).
$$
\subsection{Linearization}
Now fix a coordinate map $f\colon[n]\to[m]$.
Let $S=\Z[t_1,\ldots,t_n]$, a (commutative) polynomial ring, and let
$R=\Z[s_1,\ldots,s_m]$.  We regard $R$ as a module over $S$ via
the ring homomorphism $\wf\colon S\to R$ given by letting
$\wf(t_i)=s_{f(i)},$ for each $i$.  Let $\aug$ denote the maximal 
ideal of $S$ generated by $\set{t_1+1,t_2+1,\ldots,t_n+1}$.  We will
abuse notation and also write $\aug$ for its image $\wf(\aug)$ in $R$.
Then localize to invert the elements of $\aug$, so that
$\Z[\Z^n]\cong S_\aug$ and $\Z[\Z^m]\cong R_\aug$.

\begin{proposition}\label{prop:cw}
For any $K$, we have isomorphisms of complexes of
$\Z[\Z^n]$-modules:
\begin{eqnarray}
C^{\rm cell}_p(\Mo_K(\R,\Z))&\cong&(\Z\angl{K}^*_p\otimes S,\partial)_{\aug},
\label{eq:loc}\\
C^{\rm cell}_p(\Mo_K(\R,\Z)\times_{\Z^n}\Z^m)&\cong&(\Z\angl{K}^*_p\otimes R,
\overline{\partial})_{\aug},\text{~and}\label{eq:locR}\\
C_c^p(\Mo_K(\R,\Z)\times_{\Z^n}\Z^m)&\cong&(\Z\angl{K}_p\otimes R,
\delta)_{\aug},\label{eq:cpct}
\end{eqnarray}
for all $p\geq0$,
where the complex \eqref{eq:cpct} is the cellular cochain complex with compact 
support.  The differential $\partial$ is induced by $\partial(\vareps_i)=t_i$,
for all $i$, and $\overline{\partial}$ by $\overline{\partial}(\vareps_i)
=s_{f(i)}$.
The differential $\delta$ acts by left 
multiplication by the element $\sum_{i=1}^n e_i\otimes s_{f(i)}$.
\end{proposition}
\begin{proof}
The isomorphism \eqref{eq:loc} is obtained by localizing \eqref{eq:cell}.  
To establish \eqref{eq:cpct}, we identify $\Z[\Z^n]$ with the submodule
of $\Hom_\Z(\Z[\Z^n],\Z)$ with finite supports.  However, now $v_i$
acts (contragrediently) by multiplication by $x_i^{-1}$, for each $i$.  
Putting
this together with the dual of Proposition~\ref{prop:equivariant}(3) gives
the required isomorphism.  In fact, the right-hand side of 
\eqref{eq:cpct} is a complex of $E\otimes S_{\aug}$-modules, 
since $\delta$ commutes with multiplication by elements of $E$.
\end{proof}


Since localization is exact, it is equivalent (but more convenient) to
regard the cellular (co)chain complexes above over the polynomial ring
$S$, and we shall do so in the rest of the paper.
\begin{remark}
Our restriction to coordinate subgroups in place of arbitrary subgroups
with abelian quotient is needed for the isomorphisms of the Proposition above.
More generally, the cellular (co)chain complexes are filtered (rather
than graded) by powers
of the augmentation ideal, and the isomorphisms
\eqref{eq:loc}--\eqref{eq:cpct} are merely isomorphisms of associated graded
modules.  As Stefan Papadima observes~\cite{Pa06}, the spectral sequence
of the filtration fails to converge strongly even for very simple examples of 
non-coordinate homomorphisms.
\end{remark}

\subsection{Commutative algebra}\label{ss:commalg}
Given the setup above, it is natural to interpret group (co)homology in
terms of skew-commutative algebra.
We follow the grading conventions of \cite{eisenbook05}; in particular,
for a graded module $M$, let $M(r)_q=M_{r+q}$ for all $q$.  For modules
$M$ and $N$, the notation $\Ext^{pq}(M,N)$ refers to cohomological degree
$p$ and polynomial degree $q$.  Here, typically $q\leq-p$, so we write
$\Ext^p(M,N)_r=\Ext^{pq}(M,N)$, where $r=-p-q$.

\begin{theorem}\label{th:univ}
For any simplicial complex $K$ and for all $q\geq0$,
$$
H_q(\Mo_K(\R,\Z),\Z)\cong \Ext_E(\Z\angl{K},\Z)_q
$$
as $S$-modules.  Moreover, for all $p\geq0$,
\begin{equation}\label{eq:Hunivab}
\gr_p H_*(\Mo_K(\R,\Z),\Z)\cong\Ext_E^p(\Z\angl{K},\Z)_{\aug}
\end{equation}
where $\gr_\bfcdot$ denotes the grading associated to the filtration by
powers of the augmentation ideal of $\Z[\Z^m]$.

\end{theorem}
Before beginning the proof, we note that the (left) $S$-module structure 
on $\Ext_E(\Z\angl{K},\Z)$ here comes the identification
$S\cong\Ext_E(\Z,\Z)$ and its natural action: see \cite{sj76}.
\begin{proof}
We compute $\Ext_E(\Z\angl{K},\Z)$ by the standard injective resolution of $\Z$:
$$
0\rightarrow\Z\rightarrow I^0\rightarrow I^1\rightarrow\cdots
$$
Recall that $E$ is self-injective.  Then $\Z$ maps to the generator of 
$E_n$, and $I^q=E(n+q)\otimes_\Z S^q$, with differential 
induced $E\otimes S$-linearly by multiplication by $\sum_{i=1}^n e_i\otimes 
t_i$.  
Then
\begin{eqnarray*}
\Ext_E(\Z\angl{K},\Z) &=& H\,\Hom_E(\Z\angl{K},I^\bfcdot)\\
&\cong& H\,(\Hom_\Z(\Z\angl{K},\Z)\otimes S_\bfcdot)
\end{eqnarray*}
since $\Hom_\Z(-,\Z)=\Hom_E(-,E)(n)$.

Now we may 
compare this with Proposition~\ref{prop:cw}\eqref{eq:loc}
to obtain the cellular homology of $\Mo_K(\R,\Z)$.
\end{proof}
Now fix a simplicial complex $K$ with $n$ vertices, and a
coordinate map $f\colon[n]\to[m]$.  Let $A=A_f$ denote the kernel of the
homomorphism $\wf\colon S\to R$. 
Let $\a=\a_f$ denote the ideal
of $E$ generated in degree $1$ by functionals that vanish on $A_1$.  Let
$\ann\a$ denote its annihilator in $E$.  These ideals have the following
properties.
\begin{lemma}\label{lem:Aa}
For any $f\colon[n]\to[m]$, 
\begin{enumerate}
\item the ideal 
$A$ is generated in degree $1$, and $A_1$ is a 
free $\Z$-module of rank $n-m$.
\item 
the ideal $\a_f$ is generated by $m$ elements, for $1\leq j\leq m$:
$$
h_j := \sum_{i\colon f(i)=j}e_i,
$$
\item  the ideal $\ann\a_f$ is principal, generated by
$h_1h_2\cdots h_m$.\label{lem:Aa3}
\end{enumerate}
\end{lemma}
\begin{proof}
The first two assertions come from the definitions. 
The third follows from the fact that $\set{h_1,\ldots,h_m}$ are linearly
independent in $E_1$.
\end{proof}
In view of Lemma~\ref{lem:Aa}\eqref{lem:Aa3}, let $a_f=h_1h_2\cdots h_m$,
the generator of $\ann\a_f$.
To avoid complications, in what follows we will work over a coefficient
field $\k$.  
In order to state the next result, recall the following
definition.  The (combinatorial) Alexander dual of a simplicial complex 
$K$ on $[n]$ is a complex $K^\star$ on $[n]$.  By definition,
$$
K^\star=\set{\sigma\subseteq[n]\colon[n]-\sigma\not\in K}.
$$
\begin{theorem}\label{th:hisext}
Let $\k$ be a field, $K$ a simplicial complex, and
$f\colon[n]\to[m]$ a coordinate map.  Then, for all $q\geq0$,
\begin{eqnarray}
H_q(\Mo_K(S^1),\k[\Z^m])&\cong&\Ext_E(\k\angl{K},(a_f))_{q-n}\text{~and}
\label{eq:HisExt1}\\
H^q_c(\Mo_K(S^1),\k[\Z^m])&\cong&\Ext_E(J_{K^\star},(a_f))_{q-n},
\label{eq:HisExt2}
\end{eqnarray}
where $J_{K^\star}$ is the exterior monomial ideal associated with $K^\star$
(see \S\ref{ss:cw}), and $(a_f)$ is the principal ideal defined above.
\end{theorem}

\begin{proof}
Since $A_S$ is generated in degree $1$ and a $\k$-basis for $A_S^1$ 
is a regular sequence in $S$, its Koszul
complex is a linear, free resolution of $R$ over $S$.  That is, $R$ is
a Koszul module over $S$, so as left $E$-modules, 
$E/\a_f\cong\Ext_S(R,\k)$.  Koszul duality is an involution, so
$R\cong\Ext_E(E/\a_f,\k)$.
This is to say that $E/\a_f$ has a linear, free resolution
\begin{equation}\label{eq:res1}
\xymatrix{
0&E/\a_f\ar[l] & (E\otimes R^*,d)\ar[l]
}
\end{equation}
with a Koszul differential $\delta$ given by 
$\delta(1\otimes s_j^*)=\sum_{i\colon f(i)=j}e_i\otimes1$, 
extending $E$-linearly and by the Leibniz
rule on $R^*$.  Now apply $\Hom_E(-,E)$.
Since $\Hom_E(E/\a_f,E)$ is naturally identified with
$\ann\a_f=(a_f)$, this ideal
has an injective resolution $E\otimes R^\bfcdot$ with
differential given by multiplication by $\sum_{i=1}^n e_i\otimes s_{f(i)}$.
The proof of isomorphism \eqref{eq:HisExt1} 
concludes as in Theorem~\ref{th:univ}, using
Proposition~\ref{prop:cw}\eqref{eq:locR}.

The isomorphism~\eqref{eq:HisExt2} is analogous: using the fact that
$J_{K^\star}=\ann J_K=\Hom_E(\k\angl{K},\k)$, we see that the complex 
$(\k\angl{K}\otimes_\k R,\delta)$ of Proposition \ref{prop:cw}\eqref{eq:cpct}
actually computes $\Ext_E(J_{K^\star},(a_f))$.
\end{proof}

\setcounter{example}{0}
\begin{example}[continued]
Here, the coordinate map $f$ is an isomorphism, so $\a_f=E^{\geq1}$, and
$\ann\a\cong\k(-n)$, the socle of $E$.  Thus
Theorem~\ref{th:hisext} reduces to Theorem~\ref{th:univ} (with coefficients
in $\k$).
If $K$ is the clique complex of a graph $\Gamma$, then Theorem~\ref{th:univ}
says
$$
H_q(G_\Gamma,\Z[\Z^n])\cong\Ext_E(\Z\angl{K},\Z)_{q,\aug}
$$
for all $q\geq0$, as modules over $\Z[\Z^n]$.
In general, 
the homology of the universal abelian cover of the torus
complex is encoded in the minimal resolution
of a monomial ideal over an exterior algebra.
\end{example}
\begin{example}[continued]
In the case of the Bestvina-Brady group,
the kernel of the map $\wf:S\rightarrow\k[s_1]$ is generated by
$\set{t_i-t_j\colon 1\leq i<j\leq n}$.  
Then the ideal $\a$ is principal, generated by $a=\sum_{i=1}^n e_i$,
and $\a=\ann(\a)$.
So
$$
H_q(H_K,\k)=\Ext_E(\k\angl{K},(a))_{q,\aug}
$$
\end{example}
\section{Applications}
A main result from the work of \cite{bb97, mmv98} is that the
Bestvina-Brady group 
$H_K$ is $FP_n$ iff the flag complex $K$ is $n$-acyclic, yet $H_K$ is
finitely presented iff $K$ is simply connected.  
Then any acyclic, noncontractible
flag complex $K$ gives rise to a group which is not finitely presented,
yet has finite-dimensional (co)homology.
By regarding the homology groups as modules over the group algebra of the
abelianization, we can measure (in terms of Krull dimension) how far 
they are from having finite rank.

\subsection{Dimension calculations}
Returning to Example~\ref{ex:univab}, recall Proposition~2.1 of \cite{aah99}
provides a description of the bigraded Betti numbers of \eqref{eq:Hunivab}:
for $A=E$ or $A=S$,
let $\beta^A_{pq}(M)=\dim_\k\Tor^A_{pq}(M,\k)$, for an $A$-module $M$, and
let $P_M(t,u)=\sum_{p,q}\beta^A_{pq}(M)t^pu^q$.  Let $I=I_K$ denote the 
squarefree monomial ideal of $K$ in $S$.  Then
$$
P_{E/J}(t,u)=\sum_{p,q}\beta^S_{pq}(S/I)\frac{t^pu^q}{(1-t)^{p+q}},
$$
and
\begin{equation}\label{eq:hochster}
\beta^S_{pq}(S/I)=
\sum_{\stackrel{\I\subseteq [n]\colon
}{\abs\I=p+q}}\dim_\k 
\widetilde{H}^{q-1}(K_\I,\k),
\end{equation}
by Hochster's formula.

Then Theorem~\ref{th:univ} has the following corollaries. 
\begin{corollary}\label{cor:one}
For $q>0$, the Krull dimension of $H_q(\Mo_K(\Z,\R),\k)$ as a $\k[\Z^n]$-module
is equal to the size of the largest set $\I\subseteq[n]$
for which $\widetilde{H}^{q-1}(K_\I,\k)\neq0$.  (If there is no such set,
then $H_q(\Mo_K(\R,\Z),\k)=0$.)
\end{corollary}
\begin{proof}
It follows from \cite[Theorem 4.2, Corollary 3.8]{aah99} that
the dimension of $\Ext_E(\k\angl{K},\k)_q$ is $p+q$, where $p$ is the  
largest integer for which $\beta^S_{pq}(S/I)\neq0$.  Now use formula
\eqref{eq:hochster}.  Since the submodule of a graded $S$-module annihilated
by the ideal $\aug$ is zero, localization preserves dimension, and
our conclusion follows by Theorem~\ref{th:univ}.  
\end{proof}
We may also make use of Alexander duality.  In order to draw a parallel
with the main result of \cite{JM05}, this result is stated in terms of group
cohomology.
\begin{corollary}\label{cor:two}
Let $\Gamma$ be a graph not isomorphic to a complete graph, and $K$ its
clique complex.  Then the dimension of
the $\k[\Z^n]$-module $H^q(G_\Gamma,\k[\Z^n])$ equals the largest integer
$r$ for which there exists a simplex $\sigma\in K$ with $n-r$ vertices
satisfying $\widetilde{H}_{r-q-1}(\link_K(\sigma),\k)\neq0$.  (If there
is no such simplex, $H^q(G_\Gamma,\k[\Z^n])=0$.)
\end{corollary}
\begin{proof}
Since $\k\angl{K^\star}=E/J_{K^\star}$, it follows from the long exact
sequence for $\Ext_E$ and Theorem \ref{th:hisext} that, for all $q>0$, 
\begin{eqnarray*}
H^q(G_\Gamma,\k[\Z^n])&\cong&\Ext_E(J_{K^\star},\k)_q\\
&\cong& \widetilde{H}_{q-1}(\Mo_{K^\star}(S^1),\k[\Z^n]).
\end{eqnarray*}
Now we use the Alexander dual formulation of Hochster's formula, for
which we refer to \cite{MSbook}: if $\sigma\in K$ and $\I=[n]-\sigma$,
then for all $q$,
$$
\widetilde{H}^{q-2}(K^\star_{\I},\k)
\cong\widetilde{H}_{\abs{\I}-q-1}(\link_K(\sigma),\k).
$$
Then formula \eqref{eq:hochster} becomes
$$
\beta_{p,q-1}^S(S/I_{K^\star})=\sum_{\stackrel{\sigma\in K\colon}{
\abs{\sigma}=n-(p+q-1)}}\dim_\k \widetilde{H}_{p-2}(\link_K(\sigma),\k),
$$
from which the result follows as in Corollary~\ref{cor:one}, by letting
$r=p+q-1$.
\end{proof}
Recall that $\Z^n$ is the abelianization of $G_\Gamma$.  The Corollary
shows (by comparing with \cite{JM05}) that $H^\bfcdot(G_\Gamma,\k[G_\Gamma])$
and $H^\bfcdot(G_\Gamma,\k[G_\Gamma^{ab}])$ depend on $\Gamma$ in the
same way.

\subsection{The rank variety of $\k\angl{K}$}
Recall that if $M$ is a graded $E$-module and $a\in E^1$, then
$M$ may be regarded as a chain complex with differential given by 
multiplication by $a$.  Aramova, Avramov, and Herzog~\cite{aah99}
define an element $a\in E^1$ to be $M$-singular if the cohomology
of $(M,a)$ is nonzero, and let $V(M)$ be the variety of all
$M$-singular elements.  Say $a\in E^1$ is $M$-regular if $a\not\in V(M)$.

If $a=\sum_{i=1}^n\alpha_ie_i$, let $\supp(a)=\set{i\colon
\alpha_i\neq0}$, and consider the cohomology of $
\k\angl{K}=E/J$, regarded as a 
chain complex with differential given by multiplication by the element $a$.
We recall the following Proposition~4.3 of \cite{aah99} (with indexing
corrected).  (See also \cite[Theorem~5.5]{pasu06}.)

\begin{proposition}\label{prop:aah}
The cohomology of $(E/J,a)$ depends only on $\supp(a)$.  Let $\I=
\supp(a)$.  Then
$$
H^q(E/J,a)\cong\bigoplus_{\sigma\in K\colon\sigma\cap\I=
\emptyset}\widetilde{H}^{q-\abs{\sigma}-1}(\link_{K_\I}\sigma,
\k),
$$
where by definition
$\link_{K_\I}\sigma=\set{\tau\in K_\I\colon
\tau\cup\sigma\in K}$.
\end{proposition}
In particular, this characterizes the $E/J$-singular elements.

A sequence of elements $a_1,\ldots,a_r\in E^1$ is called $M$-regular
if $a_1$ is $M$-regular and $a_i$ is $M/(a_1,\ldots,a_{i-1})$-regular 
for each $i$, $2\leq i\leq r$.
\begin{proposition}\label{prop:finite}
For any coordinate map $f\colon[n]\to[m]$, 
$H_*(\Mo_K(S^1),\k[\Z^m])$ is finite-dimensional (over $\k$) if and only if
the elements $\set{h_1,\ldots,h_m}$ defined in Lemma~\ref{lem:Aa}(2)
form a $\k\angl{K}$-regular sequence.
\end{proposition}
\begin{proof}  
If $m=1$, $h_1=a=\sum_{i=1}^ne_i$.  Then
the ideal $(a)$ has an injective resolution
$$
\xymatrix{
0\ar[r] & (a)\ar[r] & E\ar[r]^{a} & E\ar[r]^{a} & \cdots
}
$$
so $\Ext_E(\k\angl{K},(a))$ is the cohomology of the infinite complex 
$$
\xymatrix{
E/J\ar[r]^{a} & E/J\ar[r]^{a} & \cdots
}
$$
It follows $\Ext_E^0(\k\angl{K},(a))=E/(J+(a))$, and 
$$
\Ext_E^p(\k\angl{K},E/(a))_q\cong H^q(E/J,a)
$$
for all $p\geq1$.  In particular, $\Ext_E(\k\angl{K},(a))$ 
is finite-dimensional iff $a$ is $E/J$-regular.

The general result follows from induction on $m$, using the K\"unneth formula.
\end{proof}
\setcounter{example}{1}
\begin{example}[continued]
For the Bestvina-Brady group (where $f\colon[n]\to[1]$), the only summand in 
Proposition~\ref{prop:aah} is indexed by the empty simplex, so
$H^q(E/J,a)=\widetilde{H}^{q-1}(K,\k)$.  By the argument above, 
$H^q(\Mo_K(S^1),\k[\Z])$ is finite-dimensional if and only if
$\widetilde{H}^{q-1}(K,\k)=0$, in which case 
$H^q(\Mo_K(S^1),\k[\Z])\cong E/(J+(a))^q$.  This recovers the additive
part of a result due to
Leary and Saadeto\u glu~\cite{learsaad06}.  
\end{example}


\section{Duality}\label{sec:duality}
We single out the case of Cohen-Macaulay complexes $K$ for their
particularly nice properties.  First, recall the Theorem of Eagon and 
Reiner in~\cite{ER98}.  For this, let $I_K$ denote the defining ideal of
the Stanley-Reisner ring of a complex $K$, and write $\k[K]=S/I_K$.
\begin{theorem}[\cite{ER98}]
The ideal $I_{K^\star}$ has a linear, free resolution over $S$ 
if and only if the complex $K$ is Cohen-Macaulay.
\end{theorem}

If $K$ is a Cohen-Macaulay simplicial complex of dimension $d$,
consider the Cartan complex: 
this is the complex $(\k\angl{K}\otimes S,\omega)$, 
whose differential is given by multiplication by the element
$\omega=\sum_{i=1}^n e_i\otimes u_i$.
Let $F_K=H^{d+1}(\k\angl{K}^\bfcdot\otimes S,\omega)$. 

The next Lemma is a basic consequence of Bernstein-Gelfand-Gelfand or Koszul
duality for
the module $F_K$.  We refer to \cite{eisenbook05,EFS03,ro01} for background.
\begin{lemma}\label{lem:CM}
The following are equivalent.
\begin{enumerate}
\item $K$ is Cohen-Macaulay of dimension $d$.
\item $H^p(\k\angl{K}\otimes S,\omega)=0$ for $p\neq d+1$.
\item $\Tor^S(F_K,\k)\cong\k\angl{K}^*(n-d+1)\cong J_{K^\star}(n-d+1)$ 
as graded $E$-modules.
\end{enumerate}
\end{lemma}
An interesting special case is that of $K$ a homology sphere (Gorenstein$^*$
complex).  From \cite{ER98}, it follows that:
\begin{lemma}\label{lem:gor}
$F_K\cong I_{K^\star}$ as $S$-modules.
\end{lemma}

Then BGG duality give the following reformulation
of Theorem~\ref{th:hisext}.
\begin{theorem}
For any Cohen-Macaulay
complex $K$ of dimension $d$ and  coordinate map $f\colon[n]\to[m]$,
there is an isomorphism of $S$-modules for $q\geq0$:
$$
H_q(\Mo_K(S^1),\k[\Z^m])\cong \Ext^{d+1-q}_S(F_K,R).
$$
Dually,
\begin{equation}
H^q_c(\Mo_K(S^1),\k[\Z^m])\cong \Tor^S_{d+1-q}(F_K,R).\label{eq:CMcpct}
\end{equation}
\end{theorem}
\begin{proof}
By Lemma~\ref{lem:CM}(2), the Cartan complex is a free resolution of
$F_K$ over $S$.  By applying $\Hom_S(-,R)$, we obtain the cellular
chain complex for $\Mo_K(S^1)\times_{\Z^n}\Z^m$ of 
Proposition~\ref{prop:cw}\eqref{eq:locR}, and the result follows.

The dual statement is proven in the same way.
\end{proof}
\begin{corollary}
Suppose $K$ is the clique complex of a graph $\Gamma$ and $K$ is Cohen-Macaulay
of dimension $d$. For any coordinate map $f\colon[n]\to[m]$, we have
$$
H^q(G_\Gamma,\k[\Z^m])=0\quad\text{for $q\leq m-n+d$ and $q>d+1$.}
$$
\end{corollary}
\begin{proof}
Since ${\rm pd}_S R=n-m$, $\Tor^S_i(F_K,R)=0$ for $i>n-m$.  Now
apply \eqref{eq:CMcpct}.
\end{proof}
\setcounter{example}{0}
\begin{example}[continued]
Suppose $K$ is the clique complex of $\Gamma$, and $K$ is Cohen-Macaulay
of dimension $d$.  Then the
cohomology 
with compact support of the universal abelian cover is concentrated in
dimension $d+1$; more precisely, as modules over $\k[G^{ab}_\Gamma]$,
$$
H^q(G_\Gamma,\k[G_\Gamma^{ab}])\cong\begin{cases}
(F_K)_{\aug}&\text{For $q=d+1$;}\\
0&\text{otherwise.}
\end{cases}
$$
It should be noted, however, that the vanishing of cohomology here and
in the Corollary above also follow from the paper of Brady and 
Meier~\cite{BM01}: they establish the more general result that
$G_\Gamma$ is a duality group
if and only if $K$ is a Cohen-Macaulay complex.
\end{example}
\begin{example}[continued]
Under the same hypotheses, the homology of the 
Bestvina-Brady group is given for $q\geq0$ by
$$
H_q(H_K,\k)\cong \Ext^{d+1-q}_S(I_{K^\star},\k[s]).
$$
\end{example}

\begin{ack}
The author would like to thank MSRI for its hospitality and support 
during this project.
\end{ack}

\bibliographystyle{amsplain}
\providecommand{\bysame}{\leavevmode\hbox to3em{\hrulefill}\thinspace}
\providecommand{\MR}{\relax\ifhmode\unskip\space\fi MR }
\providecommand{\MRhref}[2]{%
  \href{http://www.ams.org/mathscinet-getitem?mr=#1}{#2}
}
\providecommand{\href}[2]{#2}

\end{document}